# A Butterfly-Based Direct Integral Equation Solver Using Hierarchical LU Factorization for Analyzing Scattering from Electrically Large Conducting Objects

Han Guo[†], Yang Liu[†], Jun Hu[#], and Eric Michielssen[†]

**Abstract**

A butterfly-based direct combined-field integral equation (CFIE) solver for analyzing scattering from electrically large, perfect electrically conducting objects is presented. The proposed solver leverages the butterfly scheme to compress blocks of the hierarchical LU-factorized discretized CFIE operator and uses randomized butterfly reconstruction schemes to expedite the factorization. The memory requirements and computational cost of the direct butterfly-CFIE solver scale as $O(N\log^2 N)$ and $O(N^{1.5}\log N)$, respectively. These scaling estimates permit significant memory and CPU savings when compared to those realized by low-rank (LR) decomposition-based solvers. The efficacy and accuracy of the proposed solver are demonstrated through its application to the analysis of scattering from canonical and realistic objects involving up to 14 million unknowns.

## 1 Introduction

Electromagnetic scattering from large-scale perfect electrically conducting (PEC) objects can be analyzed using both iterative and direct surface integral equation (IE) techniques. Iterative techniques that leverage multilevel fast multipole algorithms (MLFMA) [1] or Butterfly methods [2-5] (also known as multilevel matrix decomposition algorithms) to rapidly apply discretized IE operators to trial solution vectors require $O(KN\log^\beta N)$ ($\beta=1$ or $2$) CPU and memory resources; here, $N$ is the dimension of the discretized IE operator and $K$ is the number of iterations required for convergence. The widespread adoption and success of fast iterative methods for solving real-world electromagnetic scattering problems can be attributed wholesale to their low computational costs. Iterative techniques are no panacea, however. They are ill-suited for ill-conditioned problems requiring large $K$ (e.g. scatterers supporting high-Q resonances or discretized using multi-scale/dense meshes). They also are ineffective when applied to scattering problems involving multiple excitations requiring a restart of the iterative solver for each right hand side (RHS) (e.g., calculation of monostatic radar cross section (RCS)).

[†]Department of Electrical Engineering and Computer Science, University of Michigan, Ann Arbor, MI (liuyangz@umich.edu, hanguo@umich.edu, emichiel@umich.edu).
[#]Department of Microwave Engineering, University of Electronic Science and Technology of China, Chengdu, Sichuan, China (hujun@uestc.edu.cn)

Direct methods construct a compressed representation of the inverse of the discretized IE operator and hence do not suffer (to the same degree) from the aforementioned drawbacks. Most direct methods proposed to date replace judiciously constructed blocks of the discretized IE operator and its inverse by low rank (LR) approximations [6-12]. LR compression schemes provably lead to low-complexity direct solvers for electrically small [12, 13], elongated [14, 15], quasi-planar [16], and convex [17] structures. However, for electrically large and arbitrarily shaped scatterers, the blocks of the discretized IE operators and its inverse are not LR compressible. As a result, little is known about the computational costs of LR schemes applied in this regime; experimentally their CPU and memory requirements have been found to scale as $O(N^\alpha \log^\beta N)$ ($\alpha = 2.0 \sim 3.0$, $\beta \geq 1$) and $O(N^\alpha \log N)$ ($\alpha = 1.3 \sim 2.0$), respectively.

This paper presents a low-complexity butterfly-based, direct combined field integral equation (CFIE) solver for analyzing scattering from arbitrarily shaped, electrically large 3D PEC objects. Butterfly schemes historically were developed to compress off-diagonal blocks of discretized *forward* IE operators that are not LR compressible. Here, this concept is extended, without formal proof, to inverse IE operators. This work builds on findings reported in [18, 19] that demonstrate the compressibility of discretized inverse 2D EFIE operators. Starting from a butterfly-compressed representation of the forward discretized IE operator, the proposed solver constructs a hierarchical butterfly representation of the operator's LU factors, representing all intermediate partial LU factors in terms of butterflies. The latter is achieved using a new randomized scheme to efficiently construct butterfly representations of compositions of already butterfly-compressed blocks. The CPU and memory requirements of the new solver are theoretically estimated and experimentally validated to be $O(N^{1.5} \log N)$ and $O(N \log^2 N)$, respectively. The resulting parallel direct butterfly-CFIE solver is capable of analyzing scattering from electrically large, canonical and real-life objects involving up to 14 million unknowns on a small CPU cluster.

## 2  Formulation

This section describes the proposed direct solver. Section 2.1 reviews the CFIE for PEC scatterers and its discretization. Sections 2.2 and 2.3 elucidate the butterfly scheme for compressing the discretized forward IE operator and its hierarchical LU factorization. Section 2.4 elucidates randomized butterfly reconstruction schemes for rapidly achieving the factorization.

### 2.1 Combined Field Integral Equation

Let Γ denote an arbitrarily shaped closed PEC surface residing in free space. Time harmonic electromagnetic fields $\{\boldsymbol{E}^{inc}(\boldsymbol{r}), \boldsymbol{H}^{inc}(\boldsymbol{r})\}$ impinge on Γ and induce a surface current $\mathbf{J}(\mathbf{r})$ that in turn generates scattered electromagnetic fields. Enforcing electromagnetic boundary conditions on $\Gamma$ yields the following electric and magnetic field integral equations (EFIE and MFIE):

$$\hat{\mathbf{n}} \times \hat{\mathbf{n}} \times ik\eta \int_\Gamma d\mathbf{r}' \mathbf{J}(\mathbf{r}') \cdot \left( \mathbf{I} - \frac{\nabla \nabla'}{k^2} \right) g(\mathbf{r},\mathbf{r}') = -\hat{\mathbf{n}} \times \hat{\mathbf{n}} \times \mathbf{E}^{inc}(\mathbf{r}) \quad (1)$$

$$\frac{\mathbf{J}(\mathbf{r})}{2} - \hat{\mathbf{n}} \times P.V. \int_\Gamma d\mathbf{r}' \mathbf{J}(\mathbf{r}') \times \nabla' g(\mathbf{r},\mathbf{r}') = \hat{\mathbf{n}} \times \mathbf{H}^{inc}(\mathbf{r}). \quad (2)$$

Here $\mathbf{r} \in \Gamma$, $k$ and $\eta$ denote the wavenumber and wave impedance in free space, $\hat{\mathbf{n}}$ denotes the outward unit normal to $\Gamma$, $\mathbf{I}$ is the identity dyad, $P.V.$ denotes Cauchy principal value, and $g(\mathbf{r},\mathbf{r}') = \exp(ikR/R)/(4\pi)$ with $R = |\mathbf{r} - \mathbf{r}'|$ is the free space Green's function. Both the EFIE and MFIE suffer from internal resonances. The CFIE linearly combines the EFIE and MFIE as $\text{CFIE} = \alpha \cdot \text{EFIE} + (1-\alpha)\eta \cdot \text{MFIE}$ with $0 \leq \alpha \leq 1$ and is devoid of internal resonances.

To numerically solve the CFIE, $\mathbf{J}(\mathbf{r})$ is discretized using $N$ basis functions as

$$\mathbf{J}(\mathbf{r}) = \sum_{n=1}^{N} I_n \mathbf{f}_n(\mathbf{r}). \quad (3)$$

Here $I_n$ is the expansion coefficient associated with basis function $\mathbf{f}_n(\mathbf{r})$. For simplicity, we choose the $\mathbf{f}_n(\mathbf{r})$ to be Rao-Wilton-Glisson (RWG) functions [20]. Upon Galerkin testing the CFIE, a $N \times N$ linear system of equations is obtained:

$$\mathbf{Z} \cdot \mathbf{I} = \mathbf{V}. \quad (4)$$

The $n^{th}$ entry of the solution vector $\mathbf{I}$ is $I_n$, and the $m^{th}$ entry of the excitation vector $\mathbf{V}$ is

$$V_m = \alpha \int_\Gamma d\mathbf{r} \mathbf{f}_m(\mathbf{r}) \cdot \mathbf{E}^{inc}(\mathbf{r}) + (1-\alpha)\eta \int_\Gamma d\mathbf{r} \mathbf{f}_m(\mathbf{r}) \cdot \hat{\mathbf{n}} \times \mathbf{H}^{inc}(\mathbf{r}). \quad (5)$$

The $(m,n)^{th}$ entry of the impedance matrix $\mathbf{Z}$ is

$$\begin{aligned}
\mathbf{Z}_{mn} = &-ik\alpha\eta \int_\Gamma d\mathbf{r} \int_\Gamma d\mathbf{r}' g(\mathbf{r},\mathbf{r}') [\mathbf{f}_m(\mathbf{r}) \cdot \mathbf{f}_n(\mathbf{r}')] \\
&+ \frac{i\alpha\eta}{k} \int_\Gamma d\mathbf{r} \int_\Gamma d\mathbf{r}' g(\mathbf{r},\mathbf{r}') [\nabla \cdot \mathbf{f}_m(\mathbf{r}) \nabla' \cdot \mathbf{f}_n(\mathbf{r}')] \\
&+ (1-\alpha)\eta \hat{\mathbf{n}} \cdot \int_\Gamma d\mathbf{r} \mathbf{f}_m(\mathbf{r}) \times \int_\Gamma d\mathbf{r}' \mathbf{f}_n(\mathbf{r}') \times \nabla' g(\mathbf{r},\mathbf{r}') \\
&+ \frac{(1-\alpha)\eta}{2} \int_\Gamma d\mathbf{r} \mathbf{f}_m(\mathbf{r}) \cdot \mathbf{f}_n(\mathbf{r}).
\end{aligned} \quad (6)$$

Direct solution of matrix equation (4) via Gaussian elimination or LU factorization is prohibitively expensive for large problems as these methods require $O(N^3)$ and $O(N^2)$ CPU and memory resources, respectively.

## 2.2 Butterfly Compression of the Impedance Matrix

The proposed direct solver requires a butterfly-compressed representation of $\mathbf{Z}$. The process for constructing this representation consists of two phases [2]: (i) recursively decomposing $\Gamma$ into subscatterers and (ii) compressing submatrices of $\mathbf{Z}$ representing interactions between well-separated subscatterers using butterflies.

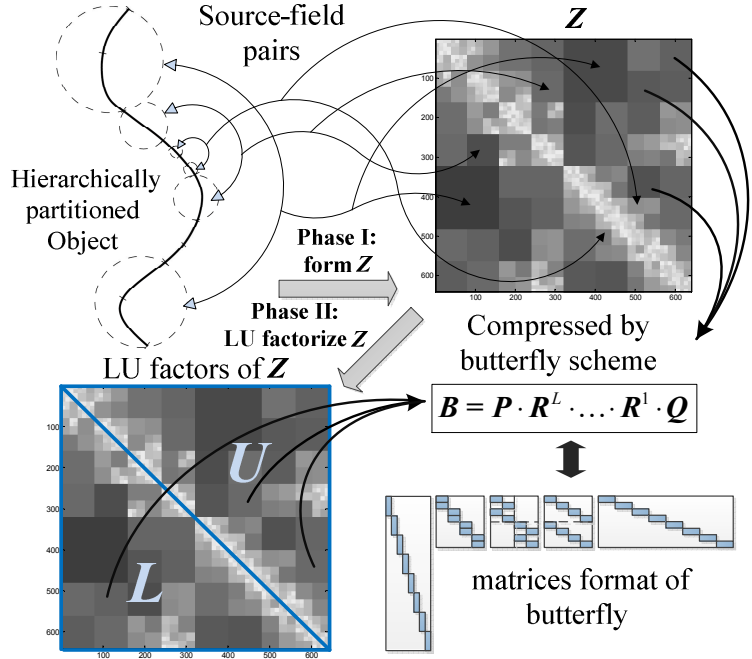

Figure 1: Matrix format of butterfly-based direct solver.

Phase (i) starts by splitting $\Gamma$ into two, roughly equal-sized level-1 subscatterers, each containing approximately $N/2$ basis functions. This splitting operation is repeated $L^h - 1$ times until the smallest subscatterers thus obtained contain $O(1)$ basis functions, resulting in a binary tree with $L^h$ levels. At level $1 \le \ell \le L^h$, there are approximately $2^\ell$ subscatterers, each containing roughly $N/2^\ell$ unknowns. Starting from level 2, two same-level subscatterers form a far-field pair if the distance between their geometric centers exceeds $2 \le \chi \le 4$ times the sum of their circumscribing radii and none of their respective ancestors constitute a far-field pair. Two level-$L^h$ subscatterers constitute a near-field pair if they do not form a far-field pair; in addition, each level-$L^h$ subscatterer forms a near-field pair with itself. This process induces a partitioning on $Z$: each level-$\ell$ far-field pair of subscatterers relates to two level-$\ell$, approximately square off-diagonal submatrices that map sources in one subscatterer to their fields observed on the other subscatterer and vice versa; additionally, each near-field pair corresponds to one or two level-$L^h$ submatrices mostly residing on or near the diagonal. The resulting decomposition is illustrated in Figure 1 assuming a simple 2D scatterer.

During Phase (ii), each far-field submatrix in $Z$ is butterfly compressed as outlined next [2, 21]. Consider a $m \times n$ level-$\ell$ far-field submatrix $Z_S^O$ with $n \approx m \approx N/2^\ell$ that models interactions between two subscatterers: a source group $S$ and an observation group $O$. This submatrix is compressed using a $L = L^h - \ell$ level butterfly. Specifically, each subscatterer is recursively subdivided using the above-described binary scheme: at level 0, there are $2^L$ source subgroups of size $n/2^L$ and one observation group of size $m$ (i.e., the group $O$ itself); at level 1, two level-0 source subgroups are paired into one level-1 source subgroup and the level-0 observation group is split into two level-1 observation subgroups; this procedure is repeated $L$ times. At level $0 \le d \le L$, there exist source subgroups $S_i^d$, $i=1,...,2^{L-d}$ and observation subgroups $O_i^d$, $i=1,...,2^d$. It follows from "degree of freedom" arguments [22] that the numerical ranks of interactions between level-$d$ subgroups are approximately constant. The maximum numerical rank $r$ for all levels is henceforth called the butterfly rank. Using this LR property, the butterfly representation $\boldsymbol{B}$ of $\boldsymbol{Z}_S^O$ consists of the product of $L+2$ sparse matrices

$$\boldsymbol{B} = \boldsymbol{P} \cdot \boldsymbol{R}^L \cdot ... \cdot \boldsymbol{R}^1 \cdot \boldsymbol{Q} \qquad (7)$$

where $\boldsymbol{P}$ and $\boldsymbol{Q}$ are block diagonal projection matrices

$$\boldsymbol{P} = \begin{pmatrix} \boldsymbol{P}_1 & & 0 \\ & \ddots & \\ 0 & & \boldsymbol{P}_{2^L} \end{pmatrix}, \quad \boldsymbol{Q} = \begin{pmatrix} \boldsymbol{Q}_1 & & 0 \\ & \ddots & \\ 0 & & \boldsymbol{Q}_{2^L} \end{pmatrix} \qquad (8)$$

with blocks of approximate dimensions $(m/2^L) \times r$ and $r \times (n/2^L)$, respectively. The interior matrices $\boldsymbol{R}^d$, $d=1,...,L$, consist of blocks of approximate dimensions $r \times 2r$ and are block diagonal following a row permutation:

$$\boldsymbol{D}^d \boldsymbol{R}^d = \begin{pmatrix} \boldsymbol{R}^{d,1} & & 0 \\ & \ddots & \\ 0 & & \boldsymbol{R}^{d,2^{L-1}} \end{pmatrix}. \qquad (9)$$

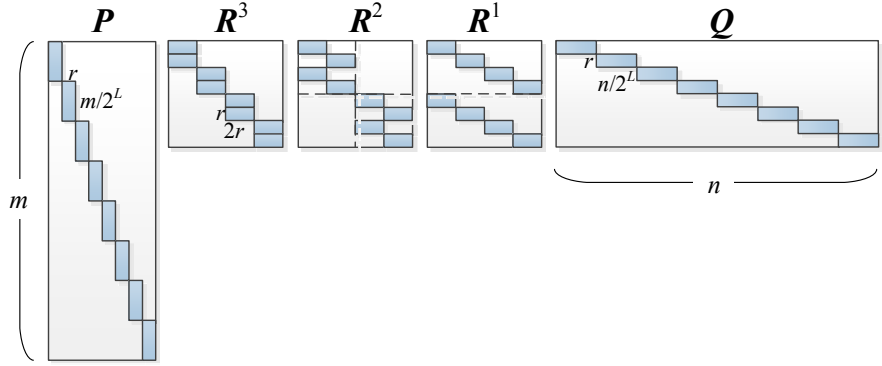

Figure 2: Algebraic structure of a 3-level butterfly

Here, $\boldsymbol{D}^d$ is the permutation matrix that renders $\boldsymbol{R}^d$ block diagonal. The structure of a 3-level butterfly is illustrated in Figure 2. In (8), the diagonal blocks $\boldsymbol{P}_i$ and $\boldsymbol{Q}_i$, $i=1,...,2^L$ are computed as [21]

$$P_i = Z_{S_{i,1}^L}^{O_i^L}, \quad Q_i = \left[ Z_{S_{1,i}^0}^{\bar{O}_i^0} \right]^\dagger Z_{S_i^0}^{\bar{O}_i^0} \tag{10}$$

and the diagonal blocks $R^{d,i}$, $i = 1,...,2^{L-1}$ in (9) are

$$R^{d,i} = \begin{pmatrix} \left[ Z_{S_{2j-1,k}^d}^{\bar{O}_{2j-1}^d} \right]^\dagger \left[ Z_{S_{j,2k-1}^{d-1}}^{\bar{O}_{2j-1}^d}, Z_{S_{j,2k}^{d-1}}^{\bar{O}_{2j-1}^d} \right] \\ \left[ Z_{S_{2j,k}^d}^{\bar{O}_{2j}^d} \right]^\dagger \left[ Z_{S_{j,2k-1}^{d-1}}^{\bar{O}_{2j}^d}, Z_{S_{j,2k}^{d-1}}^{\bar{O}_{2j}^d} \right] \end{pmatrix}. \tag{11}$$

Here, † denotes the pseudoinverse, $j = \lceil i/2^{L-d} \rceil$ and $k = \mathrm{mod}(i-1, 2^{L-d}) + 1$ where $\lceil \cdot \rceil$ and $\mathrm{mod}(\cdot,\cdot)$ denote upward rounding and modulo operations. Groups $\bar{O}_i^d \in O_i^d$, $i = 1,...,2^d$ consist of approximately $\chi_s r$ observers randomly selected from $O_i^d$ with oversampling factor $\chi_s$. Groups $S_{j,k}^d \in S_k^d$, $i = 1,...,2^d$, $k = 1,...,2^{L-d}$ consist of approximately $r$ sources that regenerate fields in level-$d$ observation subgroup $O_i^d$ due to all sources in $S_k^d$. These "skeletonized" source subgroups $S_{j,k}^d$ are constructed as follows. At level 0, $S_{1,i}^0$ is identified by performing a rank-revealing QR decomposition (RRQR) on $Z_{S_i^0}^{\bar{O}_i^0}$; at level $d > 0$, $S_{2j-1\,\mathrm{or}\,2j,k}^d$ is identified by performing a RRQR on $Z_{S_{j,2k-1}^{d-1} \cup S_{j,2k}^{d-1}}^{\bar{O}_{2j-1\,\mathrm{or}\,2j}^d}$. Note that the butterfly reduces to a LR product when $L = 0$. It is easily shown that the butterfly scheme using (10) and (11) for compressing one submatrix $Z_S^O$ requires only $O(n \log n)$ CPU and memory resources [21].

Once the impedance matrix $Z$ is hierarchically partitioned and butterfly compressed, it is hierarchically block-LU factorized as described next.

## 2.3 Hierarchical LU Factorization

Consider the following block LU factorization of $Z$:

$$Z = \begin{bmatrix} Z_{11} & Z_{12} \\ Z_{21} & Z_{22} \end{bmatrix} = \begin{bmatrix} L_{11} & \\ L_{21} & L_{22} \end{bmatrix} \begin{bmatrix} U_{11} & U_{12} \\ & U_{22} \end{bmatrix}. \tag{12}$$

In principle, the blocks featured in this factorization can be computed as follows: (i) $L_{11}$ and $U_{11}$: LU factorize $Z_{11} = L_{11} U_{11}$; (ii) $U_{12}$: solve the lower triangular subsystem $L_{11} U_{12} = Z_{12}$; (iii) $L_{21}$: solve the upper triangular subsystem $L_{21} U_{11} = Z_{21}$; (iv) $L_{22}$ and $U_{22}$: update $Z_{22} \leftarrow Z_{22} - L_{21} U_{12}$ and LU factorize the resulting $Z_{22} = L_{22} U_{22}$.

The proposed solver executes these procedures recursively until the submatrices in the LU factorized system dimension-wise match those in $Z$. More precisely, block $Z_{11}$ in step (i) and updated block $Z_{22}$ in step (iv) are decomposed similar to the original $Z$ in (12) if the corresponding blocks in $Z$ are partitioned also. Likewise, the triangular subsystems in (ii) and (iii) are further decomposed into four triangular subsystems if the corresponding $Z_{12}$ and $Z_{21}$ are partitioned in $Z$ [23]. In a similar vein, block multiplications in (iv) and in decomposition steps (ii) and (iii) are performed after decomposition of the constituent matrices if the corresponding blocks are partitioned in $Z$. Likewise, the block summations (subtractions) in (iv) and in decomposition steps (ii) and (iii) are performed after decomposition of the constituent matrices if either block is partitioned in $Z$. Following these guidelines, the hierarchical LU partitioning of $Z$ matches that of $Z$ (Figure 1). Once factorized, the inverse of the impedance matrix can be applied to excitation vectors using partitioned forward/backward substitution [12, 23].

The crux of the proposed solver lies in the experimental observation that all blocks in the above process that are not hierarchically partitioned can be butterfly compressed. The solver therefore never classically stores any ("far-field") block arising in the partial and final LU factorization of $Z$. All such blocks are processed and stored directly in butterfly-compressed format.

It is easily verified that the above-described process requires three types of "butterfly operations" that are invoked at all levels of the decomposition and consume the vast majority of the computational resources:

$$B = B_1 + B_2 \tag{13}$$

$$B = B_1 \cdot A \quad \text{or} \quad B = A \cdot B_1 \tag{14}$$

$$B = \hat{L}^{-1} \cdot B_1 \quad \text{or} \quad B = B_1 \cdot \hat{U}^{-1}. \tag{15}$$

Here, $B_1$, $B_2$ are butterfly-compressed matrices, $A$ is either a butterfly-compressed or hierarchically partitioned matrix, and $\hat{L}$ and $\hat{U}$ are hierarchically partitioned lower and upper triangular matrices consisting of butterfly-compressed blocks. The solver therefore requires a scheme for rapidly constructing butterfly-compressed versions of the matrices $B$ resulting from operations (13)-(15). In principle, such representations can be obtained using the compression scheme for the submatrices $Z_S^O$ outlined above, i.e. using (10) and (11). Unfortunately, this scheme is computationally expensive as individual entries of the $B$ matrices cannot be computed in $O(1)$ operations. That said, all $B$ s in (13)-(15) can be rapidly applied to vectors as the RHSs in (13)-(15) consist of butterfly-compressed factors/blocks (the inverse operators $\hat{L}^{-1}$ and $\hat{U}^{-1}$ are applied to vectors via partitioned forward/backward substitution [12, 23]). In what follows, two schemes for rapidly generating butterfly compressed $B$ using information gathered by multiplying $B$ and its transpose with random vectors are proposed.

**2.4 Fast Randomized Scheme for Butterfly Construction**

The two proposed randomized schemes can be regarded as generalizations of the randomized LR compression schemes recently introduced in [24, 25]. The first scheme is *iterative* in nature and permits rapid construction for butterflies with modest number of levels. The second scheme is *non-iterative* in nature and generally slower than the iterative one, but permits construction of arbitrary-level butterflies with overwhelmingly high probabilities [19, 26, 27].

*1) Iterative randomized scheme:* Consider a $L(\neq 0)$-level butterfly with dimensions $m \times n$. The proposed scheme begins by generating a $n \times n_{rnd}$ matrix $V_R$ with independent and identically distributed (i.i.d.) random entries, composed of $n_{rnd}$ column vectors. Similarly, it generates a $n_{rnd} \times m$ matrix $V_L$ with $n_{rnd}$ i.i.d. random row vectors. Let $U_R$ and $U_L^T$ denote the multiplication of the RHS and its transpose in (13)-(15) with $V_R$ and $V_L^T$, respectively:

$$U_R = B \cdot V_R \tag{16}$$
$$U_L = V_L \cdot B. \tag{17}$$

Here, the number of the random column and row vectors is chosen as
$$n_{rnd} = (L+1) \cdot r + c \tag{18}$$

with a small positive integer $c$; $r$ denotes the maximum butterfly rank (times a small oversampling factor) of all butterfly factorizations on the RHSs of (13)-(15). Furthermore, it is assumed that the projection matrices $P$ and $Q$ of the butterfly-compressed $B$ consist of blocks of dimensions $(m/2^L) \times r$ and $r \times (n/2^L)$, respectively, and that the interior matrices $R^d$, $d = 1,...L$ consist of blocks of dimensions $r \times 2r$. The proposed scheme permits rapid construction of $P$, $Q$ and $R^d$ in (13)-(15) using $U_R$, $U_L$, $V_R$ and $V_L$.

First, the projection matrices $P$ and $Q$ are computed from $U_R$ and $U_L$. Let $\bar{P}$ and $\bar{Q}$ denote initial guesses for $P$ and $Q$ obtained by filling their diagonal blocks $\bar{P}_i$ and $\bar{Q}_i$, $i = 1,...,2^L$ with i.i.d. random values. The projection matrices are computed from $\bar{U}_R = \bar{P}^T U_R$ and $\bar{U}_L = U_L \bar{Q}^T$ as
$$P_i = U_{R,i} \cdot \bar{U}_{R,i}^\dagger, \quad Q_i = \bar{U}_{L,i}^\dagger \cdot U_{L,i}, \quad i = 1,...,2^L. \tag{19}$$

Here, $U_{R,i} / \bar{U}_{R,i}$ are submatrices of $U_R / \bar{U}_R$ corresponding to the row/column indices of the block $P_i$, and $U_{L,i} / \bar{U}_{L,i}$ are submatrices of $U_L / \bar{U}_L$ corresponding to the column/row indices of block $Q_i$.

Next, the proposed scheme attempts to iteratively construct $R^d$, $d = 1,...,L$ using $\bar{U}_R$, $V_R$, $\bar{U}_L$ and $V_L$. To this end, equations (16) and (17) are rewritten as
$$\bar{U}_R = R^L \cdot .... \cdot R^1 \cdot Q \cdot V_R \tag{20}$$
$$\bar{U}_L = V_L \cdot P \cdot R^L \cdot .... \cdot R^1. \tag{21}$$

Let $U_R^d$ and $U_L^d$, $d = 0,...,L$ denote the multiplication of the partial factors in (20) and (21) with $V_R$ and $V_L$:
$$U_R^d = R^d \cdot .... \cdot R^1 \cdot U_R^0 \tag{22}$$
$$U_L^d = U_L^0 \cdot R^L \cdot .... \cdot R^{L-d+1}. \tag{23}$$

Here, $U_R^0 = Q \cdot V_R$, $U_L^0 = V_L \cdot P$, $U_R^L = \bar{U}_R$, $U_L^L = \bar{U}_L$ have already been computed. In contrast, $U_R^d$ and $U_L^d$, $d = 1,...,L-1$ are not yet known. Starting from an initial guess $R_{(0)}^d$, $d = 1,...,L$ that consists of $r \times 2r$ blocks filled with i.i.d. random values, the scheme updates $R^d$, $U_R^d$ and $U_L^d$ until convergence. Let $R_{(k)}^d$, $d = 1,...,L$ denote the updated interior matrices in the $k^{th}$ iteration. The $k^{th}$ iteration consists of approximately $L/2$ steps; each step updates interior matrices $R_{(k)}^d$ and $R_{(k)}^{L-d+1}$ for one $d = L,...,L/2+1$ as
$$U_R^{d-1} = R_{(k-1)}^{d-1} \cdot .... \cdot R_{(k-1)}^{L-d+1} \cdot R_{(k)}^{L-d} \cdot .... \cdot R_{(k)}^1 \cdot U_R^0 \tag{24}$$
$$U_R^d = R_{(k)}^d \cdot U_R^{d-1} \tag{25}$$
$$U_L^{d-1} = U_L^0 \cdot R_{(k)}^L \cdot .... \cdot R_{(k)}^d \cdot R_{(k-1)}^{d-1} \cdot .... \cdot R_{(k-1)}^{L-d+2} \tag{26}$$

$$U_L^d = U_L^{d-1} \cdot R_{(k)}^{L-d+1}. \tag{27}$$

For each $d = L,...,L/2+1$ in the $k^{\text{th}}$ iteration, the updating procedure can be summarized as follows: (i) compute $U_R^{d-1}$ using updated interior matrices $R_{(k)}^{L-d},...,R_{(k)}^1$ in (24); (ii) update $R_{(k)}^d$ by solving the linear system in (25); (iii) compute $U_L^{d-1}$ using updated interior matrices $R_{(k)}^L,...,R_{(k)}^d$ in (26); (iv) update $R_{(k)}^{L-d+1}$ by solving linear system (27). Note that $R_{(k)}^d$ and $R_{(k)}^{L-d+1}$ are block diagonal (after row permutations) (see (9)) and they can be computed by solving smaller systems (compared to those in (25) and (27)):

$$U_R^{d,i} = R_{(k)}^{d,i} \cdot U_R^{d-1,i} \tag{28}$$

$$U_L^{d,i} = U_L^{d-1,i} \cdot R_{(k)}^{L-d+1,i}. \tag{29}$$

Here $R_{(k)}^{d,i}$ and $R_{(k)}^{L-d+1,i}$, $i=1,...,2^{L-1}$ are diagonal blocks of $D^d R_{(k)}^d$ and $D^{L-d+1} R_{(k)}^{L-d+1}$, respectively. $U_R^{d,i}/U_R^{d-1,i}$ are submatrices of $U_R^d/U_R^{d-1}$ corresponding to the row/column indices of block $R_{(k)}^{d,i}$, and likewise $U_L^{d,i}/U_L^{d-1,i}$ are submatrices of $U_L^d/U_L^{d-1}$ corresponding to the column/ row indices of the block $R_{(k)}^{L-d+1}$. Specifically, the diagonal blocks $R_{(k)}^{d,i}$ and $R_{(k)}^{L-d+1,i}$ are updated as $R_{(k)}^{d,i} = U_R^{d,i} \left(U_R^{d-1,i}\right)^\dagger$ and $R_{(k)}^{L-d+1,i} = \left(U_L^{d-1,i}\right)^\dagger U_L^{d,i}$, respectively. The above-described procedure (i)-(iv) is executed $L/2$ times until all interior matrices $R_{(k)}^d$, $d = 1,...,L$ are updated. Afterwards, the scheme moves to iteration $k+1$ until the following stopping criteria is met

$$\frac{\left\| R_{(k)}^L \cdot ... \cdot R_{(k)}^1 \cdot U_R^0 - \bar{U}_R \right\| + \left\| U_L^0 \cdot R_{(k)}^L \cdot ... \cdot R_{(k)}^1 - \bar{U}_L \right\|}{\left\| \bar{U}_R \right\| + \left\| \bar{U}_L \right\|} \leq \varepsilon. \tag{30}$$

Here, $\varepsilon$ is the desired residual error. Finally, the constructed butterfly via the randomized scheme reads

$$B = P \cdot R_{(k)}^L \cdot ... \cdot R_{(k)}^1 \cdot Q. \tag{31}$$

The computational cost of the above-described scheme is dominated by that of updating $U_R^{d-1}$ and $U_L^{d-1}$, $d = L/2+1,...,L$ in each iteration via (24)/(26). These operation costs $n_{rnd} L \times O(n \log n)$ CPU time in one iteration. Therefore, the computational cost $c_n$ of the iterative randomized butterfly scheme scales as

$$c_n = k_{iter} n_{rnd} L O(n \log n) = O(n \log^3 n). \tag{32}$$

Here the convergence rate $k_{iter} = O(1)$ is not theoretically analyzed, yet we experimentally observed rapid convergence for butterflies arising in the LU factorization of CFIE operators for moderate level counts (say up to five). In the rare case where it fails to converge as may happen when butterfly level $L$ becomes large, the following non-iterative randomized scheme is used instead.

*2) Non-iterative randomized scheme:* Consider the abovementioned $L(\neq 0)$-level butterfly with dimensions $m \times n$. Just like the iterative randomized scheme, the non-iterative randomized scheme first computes the projection matrices $P$ and $Q$ using (19).

Next, the scheme computes interior matrices $R^d$, $d = 1,...,L/2$ using $U_L = V_L B$ with *structured* random vectors $V_L$. The number of random vectors is chosen as

$$n_{rnd} = r + c \tag{33}$$

Table I: Iteration counts $k_{iter}$ and measured residual errors of the iterative randomized construction scheme for butterflies with different levels $L$.

| $L$ | 1 | 2 | 3 | 4 | 5 |
|---|---|---|---|---|---|
| $k_{iter}$ | 1 | 1 | 3 | 3 | 5 |
| Error | $2.8\times10^{-9}$ | $1.3\times10^{-8}$ | $3.4\times10^{-4}$ | $1.5\times10^{-5}$ | $4.5\times10^{-4}$ |

Table II: Measured residual errors of the non-iterative randomized construction scheme for butterflies with different levels $L$.

| $L$ | 5 | 6 | 7 | 8 |
|---|---|---|---|---|
| Error | $4.2\times10^{-4}$ | $8.9\times10^{-4}$ | $1.4\times10^{-3}$ | $5.5\times10^{-3}$ |

Table III: Maximum butterfly rank in butterfly-compressed blocks at all levels of $z$ and its LU factors.

| $N$ | 6036 | 23,865 | 91,509 | 362,637 | 1,411,983 |
|---|---|---|---|---|---|
| $Z$ | 33 | 36 | 32 | 37 | 38 |
| LU of $Z$ | 33 | 36 | 39 | 46 | 50 |

with a small positive integer $c$. For each $d=1,...,L/2$ and $i=1,...,2^d$, the scheme constructs a $n_{rnd}\times m$ structured matrix $V_L$ whose columns are i.i.d. random values if they correspond to the $i^{th}$ level-$d$ observation subgroup and zero otherwise. The scheme then computes a matrix

$$V_o' = \bar{U}_L (R_{(0)}^1)^T \cdots (R_{(0)}^{d-1})^T \tag{34}$$

and a matrix $V_i' = V_o'(R_{(0)}^d)^T$, where $R_{(0)}^\ell$, $\ell=1,...,L$ is the initial guess and $\bar{U}_L = U_L \bar{Q}^T$. For each $r\times 2r$ block $R$ in $R^d$ associated with the $i^{th}$ level-$d$ observation subgroup (Figure 2), find a $n_{rnd}\times r$ submatrix $V_i$ of $V_i'$ and a $n_{rnd}\times 2r$ submatrix $V_o$ of $V_o'$ that correspond to the rows and columns of $R$. The block $R$ can be computed as $R = V_i^\dagger V_o$.

Finally, the scheme computes interior matrices $R^d$, $d=L/2+1,...,L$ using $U_R = BV_R$ with $n_{rnd}$ structured random vectors $V_R$ chosen in (33). For each $d=L,...,L/2+1$ and $i=1,...,2^{L+1-d}$, the scheme constructs a $n\times n_{rnd}$ structured random matrix $V_R$ whose rows are nonzero if they correspond to the $i^{th}$ level-$d$ source subgroup. Next, the scheme computes a matrix

$$V_o' = (R_{(0)}^{d+1})^T \cdots (R_{(0)}^L)^T \bar{U}_R \tag{35}$$

with $\bar{U}_R = \bar{P}^T U_R$. In addition, it computes a matrix $V_i' = R^{L/2} \cdots R^0 V_R$ if $d=L/2$ and $V_i' = (R_{(0)}^d)^T V_0'$ otherwise. For each $2r\times r$ block $R$ in $R^d$ associated with the $i^{th}$ level-$d$ source subgroup (Figure 2), find a $r\times n_{rnd}$ submatrix $V_i$ of $V_i'$ and a $2r\times n_{rnd}$ submatrix $V_o$ of $V_o'$ corresponding to the columns and rows of $R$. The block $R$ can be computed as $R = V_o V_i^\dagger$.

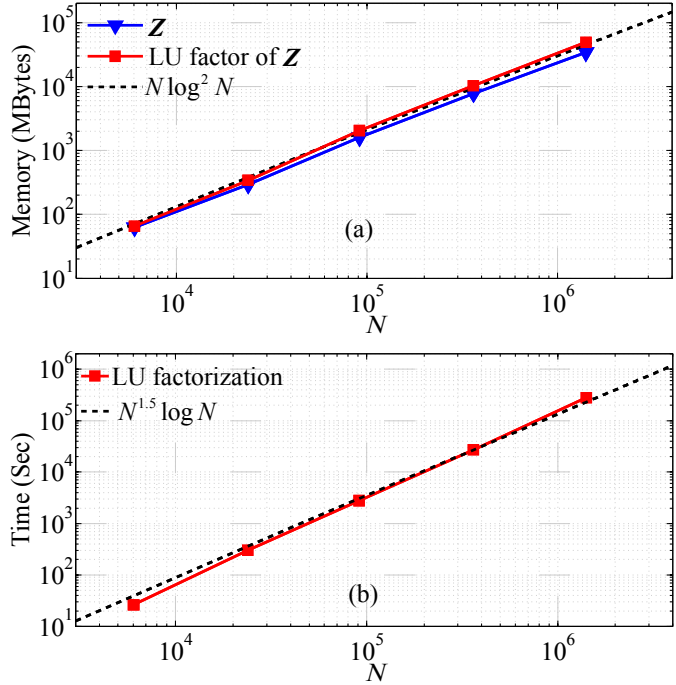

Figure 3: (a) Memory costs for storing $z$ and its LU factorization and (b) CPU times for the factorization phase using the direct butterfly-CFIE solver.

The above-described randomized scheme permits reliable construction irrespective of butterfly level $L$ given that $r$ is sufficiently large. The computational cost for the non-iterative randomized scheme $c_n$ is dominated by that of computing $U_L = V_L B$ and $U_R = B V_R$. Consequently, $c_n$ scales as

$$c_n = \sum_{d=1}^{L/2} \sum_{i=1}^{2^d} n_{rnd} O(n \log n) + \sum_{d=L}^{L/2+1} \sum_{i=1}^{2^{L+1-d}} n_{rnd} O(n \log n) = O(n^{1.5} \log n). \qquad (36)$$

Here it is assumed that multiplication of $B$ with one vector requires $O(n \log n)$ operations. Note that the scheme requires only $O(n)$ memory for storing $U_R$, $U_L$, $V_R$ and $V_L$.

## 3 Complexity Analysis

The above-described direct butterfly-CFIE solver consists of three phases: hierarchical compression of $Z$ (matrix-filling phase), hierarchical LU factorization of $Z$ (factorization phase), and application of $Z^{-1}$ to excitation vectors (solution phase). In what follows, the CPU and memory requirements of these phases are estimated.

*1) Matrix filling phase:* There are $O(2^\ell)$ level-$\ell$ far-field submatrices of approximate dimensions $(N/2^\ell) \times (N/2^\ell)$, $2 \leq \ell \leq L^h$, and the storage and CPU costs for compressing one level-$\ell$ far-field submatrix scale as $O(N/2^\ell \log(N/2^\ell))$. In addition, there are $O(N)$ level-$L^h$ near-field submatrices, each requiring $O(1)$ storage and CPU resources. Therefore, the CPU and memory requirements $C_Z$ and $M_Z$ of the matrix-filling phase scale as

$$M_Z \sim C_Z = O(N)O(1) + \sum_{\ell=2}^{L^M} O(2^\ell)O\left(\frac{N}{2^\ell}\log\frac{N}{2^\ell}\right) = O(N\log^2 N). \qquad (37)$$

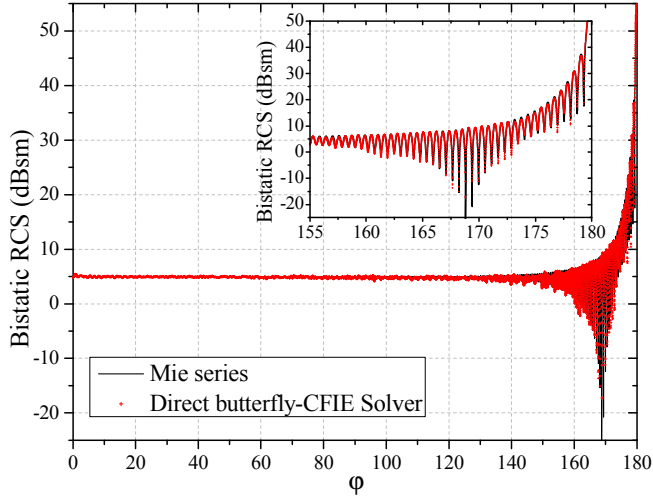

Figure 4. Bistatic RCS of the sphere at 15 GHz computed at $\theta = 90°$ and $\varphi = [0,180]°$ using the direct butterfly-CFIE solver and the Mie series.

Table IV: The technical data for the setups and solutions of the largest scattering problems considered in this paper.

|  | Sphere | Airplane |
|---|---|---|
| Maximum dimension | 2 m (100 $\lambda$) | 7.3 m (243 $\lambda$) |
| Number of unknowns $N$ | 9,380,229 | 14,179,392 |
| Number of processors | 64 | 64 |
| Memory for $Z$ | 788.2 GB | 895.1 GB |
| Memory for LU factorization | 917.4 GB | 1106.5 GB |
| Matrix filling time | 1.2 h | 1.6 h |
| Factorization time | 79.2 h | 96.6 h |
| Solving time | 43.6 s | 54.2 s |

*2) Factorization phase:* Throughout the factorization process, the number of butterfly-compressed level-$\ell$ submatrices and classically stored submatrices roughly equal those in the original impedance matrix $\mathbf{Z}$. Moreover, it is experimentally observed that $r$, representative of the sizes of the blocks in the butterfly-compressed submatrices of the LU factors, remains approximately constant during the entire factorization process. For these reasons, the memory requirement $M_{LU}$ for the factorization phase scales similarly to $M_Z$, i.e. $M_{LU} = O(N \log^2 N)$.

The CPU cost of the factorization phase is dominated by that of computing (13)-(15) while recursively carrying out procedures (i)-(iv) outlined in Section II-C. Specifically, the LU factorization of a level $\ell-1$ ($1 \leq \ell \leq L^h$) submatrix in the original and updated $\mathbf{Z}$ consists of two level-$\ell$ triangular system solves in steps (ii) and (iii), and one level-$\ell$ matrix summation (subtraction) and multiplication operation in (iv).

In steps (ii) and (iii), the level-$\ell$ off-diagonal matrix (acting as the RHS of the triangular systems) consist of at most $O(1)$ butterfly-compressed level-$d$ submatrices, $\ell+1 \leq d \leq L^h$, each corresponding to a constant number of subsystems of the form (15) and addition/multiplication operations of the form (13)/(14). Each triangular subsystem of approximate dimensions $N/2^d$, $\ell+1 \leq d \leq L^h$ can be solved in $O(c_{N/2^d})$ operations [see (15)]. Therefore, the computational cost $c_\ell^{tri}$ for solving one level-$\ell$ triangular system scales as

$$c_\ell^{tri} = \sum_{d=\ell+1}^{L^M} O(1) O(c_{N/2^d}) = O\left(\left(\frac{N}{2^\ell}\right)^{1.5} \log \frac{N}{2^\ell}\right). \tag{38}$$

Note that $c_n$ from (36) has been used in (38) as the cost of non-iterative randomized scheme dominates over that of iterative randomized scheme.

In (iv), the product of two level-$\ell$ off-diagonal matrices contains $O(1)$ butterfly-compressed level-$d$ submatrices of approximate size $N/2^d$, $\ell+1 \leq d \leq L^h$, each of which can be computed via (13) and (14) in $O(c_{N/2^d})$ operations; similarly, the updated matrix obtained by subtraction of this product from the level-$\ell$ diagonal matrix contains $O(2^{d-\ell})$ butterfly-compressed level-$d$ submatrices, $\ell+1 \leq d \leq L^h$, each of which can be computed using (13) in $O(c_{N/2^d})$ operations. Therefore, the computational costs for one level-$\ell$ matrix summation and multiplication, $c_\ell^{add}$, $c_\ell^{mul}$, scales as

$$c_\ell^{add} = \sum_{d=\ell+1}^{L^M} O(2^{d-\ell}) O(c_{N/2^d}) = O\left(\left(\frac{N}{2^\ell}\right)^{1.5} \log \frac{N}{2^\ell}\right) \tag{39}$$

$$c_\ell^{mul} = \sum_{d=\ell+1}^{L^M} O(1) O(c_{N/2^d}) = O\left(\left(\frac{N}{2^\ell}\right)^{1.5} \log \frac{N}{2^\ell}\right). \tag{40}$$

As a result, the CPU cost for the factorization phase is

$$C_{LU} = \sum_{\ell=1}^{L^M} 2 c_\ell^{tri} + c_\ell^{add} + c_\ell^{mul} = O(N^{1.5} \log N). \tag{41}$$

Here, the factor of two in front of $c_\ell^{tri}$ is due to the need to treat one upper- and one lower-triangular system.

*3) Solution phase:* In the solution phase, the computational cost $C_{RHS}$ of applying $\mathbf{Z}^{-1}$ to one RHS excitation vector via partitioned forward/backward substitution is $C_{RHS} = O(N \log N)$. No extra memory is called for.

*4) Total CPU and memory requirements:* Upon summing up the complexity estimates for the matrix filling, factorization, and solution phases, the total CPU and memory requirements $C$ and $M$ of the proposed solver scale as

$$C = C_Z + C_{LU} + C_{RHS} = O(N^{1.5} \log N) \tag{42}$$

$$M \leq M_Z + M_{LU} = O(N \log^2 N). \tag{43}$$

## 4 Numerical Results

This section presents several canonical and real-life numerical examples to demonstrate the efficiency and capabilities of the proposed direct butterfly-CFIE solver. All simulations are performed on a cluster of eight-core 2.60 GHz Intel Xeon E5-2670 processors with 4 GB memory per core. The parallelization scheme in [28] is adopted to accelerate the solver. The parallelized solver leverages a hybrid Message Passing Interface (MPI) and Open Multi-Processing (OpenMP) parallelization strategy: one MPI process is launched per processor and OpenMP utilizes eight cores on each processor.

### 4.1 Sphere

First, the performance of the randomized butterfly construction schemes is investigated. To this end, the butterfly-CFIE solver is applied to a PEC sphere of radius 1 m centered at origin. The sphere is illuminated by a *z*-polarized and *x*-propagating plane wave of 3.0 GHz. The current induced on the sphere is discretized with 565,335 RWG bases. The impedance matrix is hierarchically partitioned with 11 levels upon setting the size of the finest level groups to approximately 276 and $\chi = 2$. During the hierarchical LU factorization process, the desired residual error of the randomized scheme for constructing (13)-(15) is set to $\varepsilon = 1 \times 10^{-3}$. The typical iteration count $k_{iter}$ and the actual residual error for construction of a level-$L$ butterfly using the iterative scheme are listed in Table I. The iteration converges rapidly and the iteration count $k_{iter}$ is approximately constant. Similarly, the typical residual error for construction of a level-$L$ butterfly using the non-iterative scheme is listed in Table II.

Next, the memory requirement and the computational complexity of the direct butterfly-CFIE solver are validated. To this end, the frequency of the incident plane wave is changed from 0.3 GHz to 4.8 GHz, and the number of RWG bases is changed from 6036 to 1,411,983. The maximum butterfly ranks among all levels of $\mathbf{Z}$ and its LU factors are listed in Table III; these ranks stay approximately as constant. The memory costs for storing the impedance matrix $\mathbf{Z}$ and it factorization are plotted in Figure 3(a). These costs comply with the theoretical estimates. In addition, the CPU times of the computationally most demanding phase, v.i.z., the factorization phase, are plotted in Figure 3(b). Again, they obey the scaling estimates.

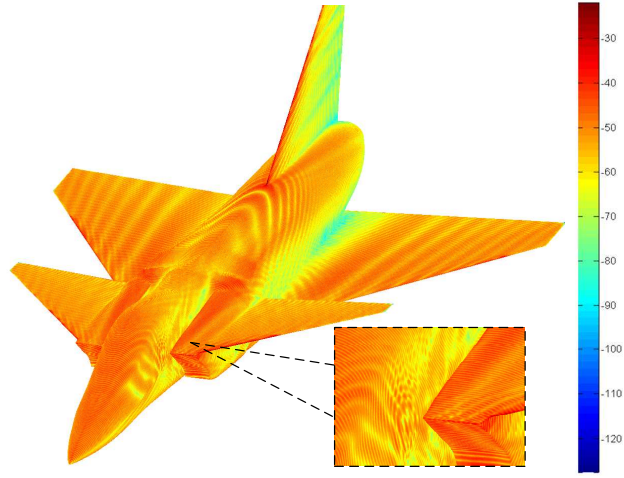

Figure 5: Current density (in dB) induced on the airplane model computed by direct butterfly-CFIE solver. The airplane is illuminated by a *z*-polarized and *x*-propagating plane wave of 10.0 GHz.

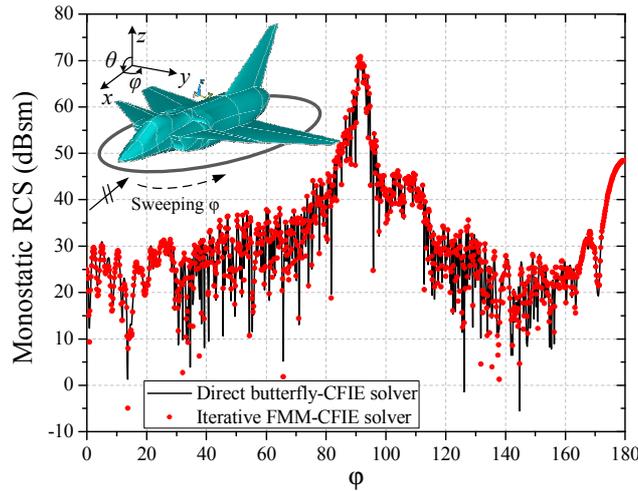

Figure 6. Monostatic RCS of the airplane at 2.5 GHz computed at $\theta = 90°$ and $\varphi = [0,180]°$ using the direct butterfly-CFIE solver and the iterative FMM-CFIE solver.

Finally, the accuracy of the direct butterfly-CFIE solver is validated by comparing bistatic radar cross section (RCS) obtained by the solver with the Mie series solutions. In this example, the sphere is illuminated by a *z*-polarized and *x*-propagating plane wave of 15 GHz. The current induced on the sphere is discretized with 9,380,229 RWG bases. The impedance matrix is hierarchically partitioned with 14 levels upon setting the size of the finest level groups to approximately 572 and $\chi = 2$. The memory costs for storing the impedance matrix and its LU factorization, and the CPU wall times for the matrix filling, factorization and solution phases are listed in Table IV. Note that $\lambda$ denotes the wavelength in Table IV. The solver requires the peak memory of 1.11 TB and total CPU time of 80.4 h using 64 processors. The bistatic RCS in directions along $\theta = 90°$ and $\varphi = [0,180]°$ are computed (Figure 4). The results agree well with the Mie series solutions.

### 4.2 Airplane Model

The capability of the direct butterfly-CFIE solver is demonstrated through its application to the scattering analysis involving an airplane model, which fits in a fictitious box of dimensions $7.30\,\text{m} \times 4.20\,\text{m} \times 1.95\,\text{m}$. The airplane is illuminated by a *z*-polarized and *x*-propagating plane wave of 10.0 GHz. The current induced on the airplane is discretized with $N = 14,179,392$ RWG bases. The impedance matrix is hierarchically partitioned with 15 levels upon setting the size of the finest level groups to approximately 432 and $\chi = 2$. The memory costs for storing the impedance matrix and its LU factorization, and the CPU wall times for the matrix filling, factorization and solution phases are listed in Table IV. The solver requires the peak memory of 1.36 TB and total CPU time of 98.2 h using 64 processors. The current induced on the airplane is shown in Figure 5.

Finally, the airplane model is illuminated by a *z*-polarized plane wave of 2.5 GHz. The plane wave is illuminating along $\theta = 90°$ and $\varphi = [0,180]°$, with a total of 10,000 incident angles. The current induced on the airplane is discretized with $N = 3,544,848$ RWG bases. The impedance matrix is hierarchically partitioned with 13 levels upon setting the size of the finest level groups to approximately 432 and $\chi = 2$. The CPU wall times using 64 processors for the matrix filling, factorization and solution phases are 14.7 min, 10.2 h and 29.5 h. Note that the averaged solving time for each RHS is only approximately 10.6 s. The peak memory cost is 268.3 GB. The monostatic RCS of the airplane is computed using the proposed solver with 10,000 incident angles and a FMM-accelerated iterative CFIE solver with 720 incident angles (Figure 6). Results are in good agreement.

## 5 Conclusions

This paper presents a butterfly-based direct CFIE solver for scattering problems involving electrically large PEC objects. The proposed solver hinges on fast randomized randomized butterfly schemes to construct a hierarchical block LU factorization of the impedance matrix. The resulting solver attains $O(N \log^2 N)$ memory and at most $O(N^{1.5} \log N)$ CPU complexities. The solver has been applied to canonical and real-life scattering problems involving millions of unknowns and many excitation vectors.

# References

[1] J. Song, C.-C. Lu, and W. C. Chew, "Multilevel fast multipole algorithm for electromagnetic scattering by large complex objects," *IEEE Trans. Antennas Propag.,* vol. 45, pp. 1488-1493, 1997.

[2] E. Michielssen and A. Boag, "A multilevel matrix decomposition algorithm for analyzing scattering from large structures," *IEEE Trans. Antennas Propag.,* vol. 44, pp. 1086-1093, 1996.

[3] E. Candes, L. Demanet, and L. Ying, "A fast butterfly algorithm for the computation of Fourier integral operators," *Multiscale Model. Sim.,* vol. 7, pp. 1727-1750, 2009.

[4] M. Tygert, "Fast algorithms for spherical harmonic expansions, III," *J. Comput. Phys.,* vol. 229, pp. 6181-6192, 2010.

[5] Y. Li and H. Yang, "Interpolative butterfly factorization," *arXiv preprint arXiv:1605.03616,* 2016.

[6] J.-G. Wei, Z. Peng, and J.-F. Lee, "A fast direct matrix solver for surface integral equation methods for electromagnetic wave scattering from non-penetrable targets," *Radio Sci.,* vol. 47, 2012.

[7] L. Greengard, D. Gueyffier, P.-G. Martinsson, and V. Rokhlin, "Fast direct solvers for integral equations in complex three-dimensional domains," *Acta Numerica,* vol. 18, pp. 243-275, 2009.

[8] A. Heldring, J. M. Rius, J. M. Tamayo, J. Parron, and E. Ubeda, "Multiscale compressed block decomposition for fast direct solution of method of moments linear system," *IEEE Trans. Antennas Propag.,* vol. 59, pp. 526-536, 2011.

[9] J. Shaeffer, "Direct solve of electrically large integral equations for problem sizes to 1 M unknowns," *IEEE Trans. Antennas Propag.,* vol. 56, pp. 2306-2313, 2008.

[10] H. Guo, J. Hu, H. Shao, and Z. Nie, "Hierarchical matrices method and its application in electromagnetic integral equations," *Int. J. Antennas Propag.,* vol. 2012, 2012.

[11] W. Chai and D. Jiao, "An-matrix-based integral-equation solver of reduced complexity and controlled accuracy for solving electrodynamic problems," *IEEE Trans. Antennas Propag.,* vol. 57, pp. 3147-3159, 2009.

[12] M. Bebendorf, "Hierarchical LU decomposition-based preconditioners for BEM," *Computing,* vol. 74, pp. 225-247, 2005.

[13] E. Corona, P.-G. Martinsson, and D. Zorin, "An $O(N)$ direct solver for integral equations on the plane," *Appl. Comput. Harmon. Anal.,* vol. 38, pp. 284-317, 2015.

[14] P.-G. Martinsson and V. Rokhlin, "A fast direct solver for scattering problems involving elongated structures," *J. Comput. Phys.,* vol. 221, pp. 288-302, 2007.

[15] E. Michielssen, A. Boag, and W. Chew, "Scattering from elongated objects: direct solution in $O(N\log^2 N)$ operations," in *IEE P-Microw. Anten. P.,* 1996, pp. 277-283.

[16] E. Winebrand and A. Boag, "A multilevel fast direct solver for EM scattering from quasi-planar objects," in *Proc. Int. Conf. Electromagn. Adv. Appl.,* 2009, pp. 640-643.

[17] Y. Brick, V. Lomakin, and A. Boag, "Fast direct solver for essentially convex scatterers using multilevel non-uniform grids," *IEEE Trans. Antennas Propag.,* vol. 62, pp. 4314-4324, 2014.

[18] H. Guo, J. Hu, and E. Michielssen, "On MLMDA/butterfly compressibility of inverse integral operators," *IEEE Antenn. Wireless Propag. Lett.,* vol. 12, pp. 31-34, 2013.

[19] Y. Liu, H. Guo, and E. Michielssen, "A HSS matrix-inspired butterfly-based direct solver for analyzing scattering from two-dimensional objects," *IEEE Antenn. Wireless Propag. Lett.,* vol. PP, pp. 1-1, 2016.

[20] S. M. Rao, D. R. Wilton, and A. W. Glisson, "Electromagnetic scattering by surfaces of arbitrary shape," *IEEE Trans. Antennas Propag.,* vol. 30, pp. 409-418, 1982.

[21] E. Michielssen and A. Boag, "Multilevel evaluation of electromagnetic fields for the rapid solution of scattering problems," *Microw. Opt. Techn. Let.,* vol. 7, pp. 790-795, 1994.

[22] O. M. Bucci and G. Franceschetti, "On the degrees of freedom of scattered fields," *IEEE Trans. Antennas Propag.,* vol. 37, pp. 918-926, 1989.

[23] S. Borm, L. Grasedyck, and W. Hackbusch, "Hierarchical matrices," *Lecture notes,* vol. 21, p. 2003, 2003.

[24] E. Liberty, F. Woolfe, P.-G. Martinsson, V. Rokhlin, and M. Tygert, "Randomized algorithms for the low-rank approximation of matrices," *Proc. Natl. Acad. Sci.,* vol. 104, pp. 20167-20172, 2007.

[25] P.-G. Martinsson, V. Rokhlin, and M. Tygert, "A randomized algorithm for the decomposition of matrices," *Appl. Comput. Harmon. Anal.,* vol. 30, pp. 47-68, 2011.

[26] Y. Li, H. Yang, E. R. Martin, K. L. Ho, and L. Ying, "Butterfly factorization," *Multiscale Model. Sim.,* vol. 13, pp. 714-732, 2015.

[27] Y. Liu, H. Guo, and E. Michielssen, "A new butterfly reconstruction method for MLMDA-based direct integral equation solvers," in *Proc. IEEE Int. Symp. AP-S/URSI,* 2016.

[28] H. Guo, J. Hu, and Z. Nie, "An MPI-OpenMP hybrid parallel-LU direct solver for electromagnetic integral equations," *Int. J. Antennas Propag.,* 2015.